\input amstex
\documentstyle{amsppt}
\magnification 1200
\NoRunningHeads
\NoBlackBoxes
\document

\def\ell{{\text{ell}}}

\def\h1{\hat{\bold 1}}

\def\hh{\hbar}

\def\Ua{U_q(\tilde\g)}
\def\U2{{\Ua}_2}
\def\g{\frak g}
\def\n{\frak n}

\def\Z{\Bbb Z}
\def\C{\Bbb C}
\def\d{\partial}

\def\<{\langle}
\def\>{\rangle}
\def\o{\otimes}

\def\End{\text{End}}

\def\h{{\frak h}}
\topmatter
\title Whittaker functions on quantum groups and q-deformed Toda operators
\endtitle

\author {\rm {\bf Pavel Etingof} \linebreak
\vskip .1in
Department of Mathematics\linebreak
Harvard University\linebreak 
Cambridge, MA 02138, USA\linebreak
e-mail: etingof\@math.harvard.edu}
\endauthor
\endtopmatter

\centerline{Dedicated to Dmitry Borisovich Fuchs on the occasion of his 
sixtieth birthday}

\head Introduction\endhead

Let $G$ be a simply connected simple Lie group over $\C$.
Let $N_\pm$ be the positive and the negative maximal unipotent 
subgroups, and $H$ the maximal torus, corresponding to some 
polarization of $G$. 
Let $G_0=N_-HN_+$ be the big Bruhat cell. 

Let 
$\chi_\pm:N_\pm\to \C^*$ be holomorphic nondegenerate 
characters (i.e. they don't
vanish on simple roots). A Whittaker function on $G_0$ with 
characters $\chi_+,\chi_-$ is any holomorphic function 
$\phi$ on $G_0$ such that $\phi(a_-a_0a_+)=
\chi_-(a_-)\chi_+(a_+)\phi(a_0)$, where $a_\pm\in N_\pm$,
$a_0\in G_0$. Thus, a Whittaker function is completely determined 
by its values on the maximal torus $H$. 

In the late 1970's it was observed (by Kazhdan and Kostant)
that the restriction of the 
Laplace operator on $G$ to Whittaker functions is the quantum Toda 
Hamiltonian. This allows one to easily prove 
Kostant's integrability theorem
for the quantum Toda 
system: quantum integrals are restrictions 
to Whittaker functions of higher Casimirs of $G$.  
 
This procedure can be generalized to the case when the group $G$ is replaced 
with the corresponding affine Kac-Moody group $\hat G$. In this case, 
one should consider Whittaker functions ``of critical level'', i.e. 
functions satisfying the equation $\phi(xz)=\phi(x)z^{-h^\vee}$, where $z$ is a
central element of $\hat G$ and $h^\vee$ the dual Coxeter number of $G$. 
The restriction of the Laplace operator to Whittaker functions
is then the quantum affine Toda Hamiltonian. 
As a result, one gets a proof of the integrability of the quantum affine 
Toda system:  
quantum integrals are restrictions 
to Whittaker functions of higher Casimirs of $\hat G$ at the critical level
defined in \cite{FF,GW}. 
 
The goal of this paper is to generalize this theory 
(both for $G$ and $\hat G$) to the case 
of quantum groups, and give the corresponding proofs of the quantum 
integrability of the q-difference analogs of the quantum Toda systems
for $G$ and $\hat G$. 

The main problem with this generalization is that the algebra 
$U_q(\frak n_+)$, which is the q-analogue 
of $N_+$, does not have a homomorphism to $\C$ which is a deformation
of $\chi_+$. This problem can be dealt with as follows.
 
Consider the universal character $\eta_+:\n_+\to \C^r$ (where 
$r$ is the rank of $G$), so that individual 
nondegenerate characters $\chi_+$ are obtained 
by the formula $\chi_+=\xi\circ \eta_+$, where $\xi:\C^r\to\C$ is a linear map 
whose coordinates are nonzero. (Here we abuse notation 
by denoting the Lie algebra homomorphism corresponding to the group 
homomorphism $\chi_+$ also by $\chi_+$).
We can regard $\eta_+$ as a homomorphism 
from $U(\n_+)$ to the polynomial ring $\C[x_1,...,x_r]$. 
The well-known but crucial observation is that $\eta_+$ 
admits a quantum deformation $\eta_+^q$, which is a map 
from $U_q(\n_+)$ to a certain quantum 
polynomial algebra. This allows one to generalize 
the definition of a Whittaker function to the q-case, 
after which it is more or less straightforward to generalize the
results about Toda systems. 

{\bf Remark.} A different (but closely 
related) method of dealing with the problem of absence of 
characters is to multiply the generators of $U_q(\n_+)$ by elements 
of the Cartan subgroup in such a way that the algebra generated by the 
obtained elements 
has nondegenerate characters. This beautiful idea was introduced in the recent 
paper \cite{S}, which appeared while our paper was being prepared.  
This approach, as indicated in \cite{S}, can also be used to produce 
q-deformations of quantum Toda systems.

The paper is structured as follows. In Section 1 we recall
the Kazhdan-Kostant construction for finite dimensional classical groups. 
In Section 2 we generalize this construction to the affine case. 
In Section 3 we generalize the construction to finite-dimensional 
quantum groups using the method described above. In Section 4 we give the 
generalization to affine quantum groups. In Section 5 we explain how to 
compute explicitly the Toda Hamiltonians for quantum and 
affine quantum groups, and do so for type A. In Section 6 we discuss 
the connection of our q-deformed Toda systems with the ordinary Toda systems
and their known q-deformations. In Section 7 we discuss the 
relation between the Toda systems and the Calogero-Moser, Macdonald, 
and Ruijsenaars systems.

{\bf Acknowledgments} The author thanks David Kazhdan for posing the problem
and Bertram Kostant for discussions. The author is also grateful 
to Victor Ginzburg for discussions and references. 

\head 1. Toda systems\endhead

We start with recalling the Kazhdan-Kostant construction. 

To avoid difficulties with generalizing the results to quantum groups, 
we will take an algebraic approach and 
use formal groups rather than ordinary complex Lie groups. 

Let $\g$ be a simple complex Lie algebra, and $\h$ a Cartan subalgebra in $\g$.
We fix an invariant inner product on $\g$ such that 
long roots have squared norm $2$. 

\proclaim{Definition} The quantum Toda Hamiltonian  
corresponding to the Lie algebra $\g$ is the
following differential operator 
on $\h$:
$$
M=-\frac{1}{2}\Delta+\sum_{i=1}^re^{\alpha_i(h)},\tag 1.1
$$
where 
$\Delta$ is the Laplacian corresponding to the invariant form, and 
$\alpha_i$ are simple positive roots of $\g$. 
\endproclaim

\proclaim{Theorem 1.1} (\cite{K}) The operator $M$ defines a completely 
integrable quantum system. More precisely, there exist differential operators 
$M_1,...,M_r$ on $\h$ ($r=rank(\g)$) such that 

(a) $[M_i,M_j]=[M_i,M]=0$;

(b) The 
symbols of $M_i$ are Weyl group invariant elements of $S\h$ which generate 
$(S\h)^W$.
\endproclaim

The proof of Theorem 1.1 occupies the rest of Section 1. 

Choose a polarization $\g=\n_+\oplus\h\oplus\n_-$. 

Let $C\in U(\g)$ be the quadratic 
Casimir corresponding to the invariant form on $\g$. We have 
$$
C=\sum_{i=1}^ry_i^2+2\sum_{\alpha>0}f_\alpha e_\alpha+2\rho,\tag 1.2
$$
where $y_i$ is an orthonormal basis of $\h$, 
$e_\alpha,f_\alpha$ are root elements, and $\rho$ is the element of $\h$ 
such that $\alpha(\rho)=1$ for all simple roots $\alpha$. 

Let $A_\g=U(\g)^*$ be the algebra of functions on the formal group $G$
associated to $\g$. 

For any $X\in U(\g)$, let $R(X),L(X)\in \End_\C(A_\g)$ be the  
right-invariant and left-invariant
differential operators on $G$ corresponding to $X$. Thus we get 
$$
R(C)=L(C)=\sum_{i=1}^rL(y_i)^2+2\sum_{\alpha>0}L(f_\alpha) L(e_\alpha)+
2L(\rho).\tag 1.3
$$

Let $A_\h=U(\h)^*$ be the algebra of functions on the formal group 
$H$ associated to $\h$. Let $\phi\to \phi|_\h$ denote the map 
$A_\g\to A_\h$ dual to the embedding $U(\h)\to U(\g)$. 
Let $\d_h$ be the derivation 
of the algebra $A_\h$ corresponding to $h\in \h$. 
In particular, let $\d_{y_i}=\d_i$.  

Let $\chi_\pm:U(\n_\pm)\to \C$ be defined by 
$\chi_+(e_i)=\chi_+^i,\chi_-(f_i)=\chi_-^i$
($e_i,f_i$ are simple root vectors). 

\proclaim{Definition} A Whittaker function with characters 
$\chi_+,\chi_-$ is a function $\phi\in A_\g$  
such that $L(X_+)\phi=\chi_+(X_+)\phi$, 
$R(X_-)\phi=\chi_-(X_-)\phi$, where $X_\pm\in U(\n_\pm)$.
\endproclaim

By Frobenius reciprocity, the space of Whittaker functions 
can be naturally identified with $A_\h$, via $\phi\to \phi|_\h$. 

Let $\phi\in A_\g$ be a Whittaker function
with characters $\chi_+$, $\chi_-$. 

\proclaim{Proposition 1.2} One has 
$$
(L(C)\phi)|_\h=(\sum_{i=1}^r \d_i^2+
2\d_{\rho}+2\sum_{i=1}^r \chi_+^i\chi_-^ie^{\alpha_i(h)})\phi|_\h
\tag 1.4
$$
(where $h\in \h$).
\endproclaim

\demo{Proof} 
Applying (1.3) to $\phi$, we get
$$
(L(C)\phi)|_\h=(\sum_{i=1}^r \d_i^2+
2\d_{\rho})\phi|_\h+(2\sum_{i=1}^r \chi_+^iL(f_i)\phi)|_\h\tag 1.5 
$$
So it remains to compute $(L(f_i)\phi)|_\h$. We have 
$L(f_i)(e^h)=Ad(e^h)^{-1}(R(f_i)(e^h))=e^{\alpha_i(h)}R(f_i)(e^h)$.
Therefore, we get (1.4).
$\square$\enddemo

Let $D'$ be the differential operator on the RHS of (1.4). 
It is easy to see that the operator 
$D=e^{(\rho,h)}D'e^{-(\rho,h)}$ has the form 
$$
D=\sum_{i=1}^r \d_i^2
+\sum_{i=1}^r 2\beta_ie^{\alpha_i(h)}+(\rho,\rho), 
\beta_i=\chi_+^i\chi_-^i.\tag 1.6
$$
Thus, if $\beta_i=-1$ then $M:=-\frac{1}{2}(D-(\rho,\rho))$ 
is the quantum Toda Hamiltonian (1.1). 

Now let us find quantum integrals of $M$. 

\proclaim{Proposition 1.3} 
(i) For any element $Y\in U(\g)$ there exists 
a unique differential operator 
$D_Y:A_\h\to A_\h$ on $H$ such that for any Whittaker function 
$\phi\in A_\g$ one has $(L(Y)\phi)|_\h=D_Y\phi|_\h$. 

(ii) If $Y$ and $Y'$ are central in $U(\g)$ then $D_{YY'}=D_YD_{Y'}$.
\endproclaim

\demo{Proof} (i) The proof is a straightforward generalization 
of the argument in the proof of Proposition 1.2. 

(ii) The statement follows from the fact that if $Y$ is central 
then $L(Y)$ maps Whittaker functions 
to themselves. 
$\square$\enddemo

Now let $Y_1,Y_2,...,Y_r$ be a system 
of generators of $U(\g)^\g$, and $M_i=e^{(\rho,h)}D_{Y_i}e^{-(\rho,h)}$.
Then $M_i$  
satisfy the conditions of Theorem 1.1, Q.E.D.  

\head 2. Affine Toda systems\endhead

In this section we generalize the results of Section 1 to 
the affine case. 

\proclaim{Definition} The affine quantum Toda Hamiltonian  
corresponding to the Lie algebra $\g$ is the
following differential operator 
on $\h$:
$$
M^K=-\frac{1}{2}\Delta+\sum_{i=1}^re^{\alpha_i(h)}+Ke^{-\theta(h)},\tag 2.1
$$
where 
$\Delta$, $\alpha_i$ are as above, 
$\theta$ is the maximal root of $\g$, and $K$ is a nonzero complex number. 
\endproclaim

{\bf Remark.} 
Note that $K$ is an essential parameter and cannot be removed by a simple 
change of variables. Also observe that $M^0=M$, 
where $M$ is as in the previous section. 

The following
generalization of Theorem 1.1 
was proved by Olshanetsky and Perelomov \cite{OP} for Lie 
algebras of type A and follows from the results of Cherednik \cite{Ch1}
(see Section 7) in the general case. 

\proclaim{Theorem 2.1} The operator $M^K$ defines a completely 
integrable quantum system. More precisely, there exist differential operators 
$M_1^K,...,M_r^K$ on $\h$ ($r=rank(\g)$) such that 

(a) $[M^K_i,M^K_j]=[M^K_i,M^K]=0$;

(b) The 
symbols of $M^K_i$ are Weyl group invariant elements of $S\h$ which generate 
$(S\h)^W$.
\endproclaim

The proof of this theorem is a generalization of the proof of Theorem 1.1, 
and is given below. The main modification in the proof, compared to
Section 1, is that now instead of the Lie algebra $\g$ one should consider the 
affine Lie algebra $\hat\g$ associated to $\g$. 

Recall \cite{Kac} that the affine Lie algebra $\hat\g$ has the form 
$L\g\oplus \C c$, where $L\g$ is the Lie algebra of Laurent polynomials 
of a variable $t$ with values in $\g$, and $c$ is a central element. 
The commutator is defined by $[f(t),g(t)]=[fg](t)+\text{Res}_{t=0}(df,g)c$. 
 
Let $\hat\n_+$ be the Lie subalgebra of $L\g$
of elements $g(t)$ regular at $0$ and such that $g(0)\in \n_+$. 
Let $\hat\n_-$ be the Lie subalgebra of $L\g$
of elements $g(t)$ regular at infinity and 
such that $g(\infty)\in \n_-$. 
Let $\hat \h=\h\oplus\C c$. Thus we have 
$\hat\g=\hat\n_+\oplus\hat\h\oplus\hat\n_-$. 

Let $A_{\hat\g}=U(\hat\g)^*$ be the algebra of functions 
on the formal group $\hat G$ associated to $\hat\g$. 
For any $X\in U(\hat\g)$, let $R(X),L(X)$ be the left- and right-invariant 
differential operators on $\hat G$ corresponding to $X$
(i.e. $R(X),L(X)$ are endomorphisms of the space $A_{\hat\g}$). 

For any complex number $k$, let $A^k_{\hat\g}\subset A_{\hat\g}$ 
be the space of functions 
satisfying the equation $L(c)f=kf$. 
Then it is clear that for any $Y\in U(\hat\g)$, 
$L(Y)$ preserves $A^k_{\hat\g}$, and $L(c)|_{A^k_{\hat\g}}=k\cdot Id$. 
Thus $L$ descends to a map 
$L:U(\hat\g)/(c-k)\to \End_\C(A^k_{\hat\g})$. 

 The value $k=-h^\vee$, where $h^\vee$ is the dual Coxeter number of $\g$
(the so called critical level) is especially important for us,
 because of the following theorem. 

Let $\hat U(\hat\g)$ 
be the completion of $U(\hat\g)$ 
acting in the category of $\hat\g$-modules
which are locally nilpotent under the action of
simple root vectors $e_i$, $i=0,...,r$
(i.e. for any vector $v$ there exists $N(v)>0$ such that 
$e_{i_1}...e_{i_{N(v)}}v=0$ for any $i_1,...,i_{N(v)}$).

\proclaim{Theorem 2.2} (\cite{FF,GW}) 
The algebra $\hat U(\hat\g)/(c+h^\vee)$ contains 
algebraically independent central elements $Y_1,...,Y_r$, 
of degree $0$ with respect to the t-grading, such that:

(a) $Y_1=\hat C:=C+2\sum_{m\ge 0}\sum_{a\in B}(a\o t^{-m})(a\o t^m)$, 
where $B$ is an orthonormal basis of $\g$; 

(b) $Y_i$ are of the form $Y_i=Y_i^0+Y_i^+$, where 
$Y_i^0$ is a Weyl group invariant element of $S\h$, 
and $Y_i^+$ is a sum of monomials 
of degree $\le \text{deg}(Y_i^0)$, which belong to 
$U(\hat\g)\hat\n_+$. 
\endproclaim

 Now we define Whittaker functions, in the same way as before. 
Let $\chi_\pm:U(\hat\n_\pm)\to \C$ be defined by 
$\chi_+(e_i)=\chi_+^i,\chi_-(f_i)=\chi_-^i$
($e_i,f_i$ are simple root vectors, $i=0,...,r$). 

\proclaim{Definition} A Whittaker function with characters 
$\chi_+,\chi_-$ is a function $\phi\in A^{-h^\vee}_{\hat\g}$  
such that $L(X_+)\phi=\chi_+(X_+)\phi$, 
$R(X_-)\phi=\chi_-(X_-)\phi$, where $X_\pm\in U(\hat\n_\pm)$.
\endproclaim

By Frobenius reciprocity, the space of Whittaker functions 
can be naturally identified with $A_\h$, via $\phi\to \phi|_\h$. 

Let $\phi\in A^{-h^\vee}_{\hat\g}$ be a Whittaker function
with characters $\chi_+$, $\chi_-$. 

\proclaim{Proposition 2.3} The series 
$((L(\hat C)\phi)|_\h$ is finite. Moreover, one has 
$$
(L(\hat C)\phi)|_\h=(\sum_{i=1}^r \d_i^2+
2\d_{\rho}+2\sum_{i=0}^r \chi_+^i\chi_-^ie^{\alpha_i(h)})\phi|_\h
\tag 2.2
$$
(where $h\in \h$).
\endproclaim

(The difference between the right hand sides of (1.4) and (2.2) is that 
in (2.2) the second summation includes $i=0$.)

\demo{Proof} The proof is completely analogous to the proof of 
Proposition 1.2. 
$\square$\enddemo

Let $\hat D'$ be the differential operator on the RHS of (2.2). 
It is easy to see that the operator 
$\hat D=e^{(\rho,h)}\hat D'e^{-(\rho,h)}$ has the form 
$$
\hat D=\sum_{i=1}^r \d_i^2
+\sum_{i=0}^r 2\beta_ie^{\alpha_i(h)}+(\rho,\rho), 
\beta_i=\chi_+^i\chi_-^i.\tag 2.3
$$
Thus, if $\beta_i=-1$ for $i\ge 1$, and $\beta_0=-K$ 
then $M^K=-\frac{1}{2}(\hat D-(\rho,\rho))$ 
is the affine quantum Toda Hamiltonian (2.1). 

Now let us find quantum integrals of $M^K$. 

\proclaim{Proposition 2.4} 
(i) For any $i=1,...,r$ the series $(L(Y_i)\phi)|_\h$ is finite
for any Whittaker function $\phi$. Moreover, 
there exists  
a unique differential operator 
$D_{Y_i}:A_\h\to A_\h$ on $H$ such that for any Whittaker function 
$\phi$ one has $(L(Y_i)\phi)|_\h=D_{Y_i}\phi|_\h$. 

(ii) $[D_{Y_i},D_{Y_j}]=0$. 
\endproclaim

\demo{Proof} The finiteness of the series follows 
from the fact that terms containing non-simple roots 
vanish. The rest of the proof is analogous to the proof of Proposition 1.3. 
$\square$\enddemo

Now let $M_i^K=e^{(\rho,h)}D_{Y_i}e^{-(\rho,h)}$.
Then $M_i^K$  
satisfy the conditions of Theorem 2.1, Q.E.D.

\head{3. q-deformed Toda systems}\endhead

In this section we will generalize the constructions of Section 1 to the case 
when the classical group $G$ is replaced with the corresponding quantum
group. 

Recall that $H$ denotes the (formal) maximal torus in $G$, and 
$A_\h$ is the algebra of regular functions on $H$. 

\proclaim{Definition} A
difference operator on $H$ is an operator on the space 
$A_\h[[\hh]]$ (where $\hh$ is a formal parameter) 
which has the form $\sum f_iT_{\beta_i}$ (a finite sum), where 
$f_i\in A_\h[[\hh]]$, and $T_\beta f(x)=f(xe^{\hh \beta})$, 
$\beta\in \h$. 
\endproclaim

Difference operators form an algebra, which we'll call $D_q(H)$. 
For any algebra $B$, by a $B$-valued difference operator we will mean 
an element of $D_q(H)\o B$ (the algebraic, i.e. uncompleted tensor product). 

Below, using 
the method of Section 1, we will produce a commuting family of 
$r$ algebraically independent scalar valued difference operators on $H$, 
which is a deformation of the Toda system.  

Let $U_\hh(\g)$ be the Drinfeld-Jimbo quantum 
universal enveloping algebra corresponding 
to the Lie algebra $\g$ 
(with the formal quantization parameter). 
It has generators $e_i$, $f_i$, $h_i$ 
and the standard relations \cite{CP}.
Let $U_\hh(\n_+)$ be the subalgebra of $U_\hh(\g)$ generated by $e_i$, 
and $U_\hh(\n_-)$ be the subalgebra generated by $f_i$, respectively. 

Let $A_{U_\hh(\g)}=U_\hh(\g)^*$ be the 
space of linear functions on $U_\hh(\g)$, i.e. the dual quantum formal group. 
The space $A_{U_\hh(\g)}$ is the quantum analog of $A_\g$.  

The obvious two-sided action of $U_\hh(\g)$ on itself 
induces a two-sided action of $U_\hh(\g)$ on $A_{U_\hh(\g)}$, via
$((a,b)\circ \phi)(x)=\phi(axb)$ (here the action of $a$ is a right 
action and the action of $b$ is a left action). Define 
the maps $L,R:U_\hh(\g)\to 
\End(A_{U_\hh(\g)})$ by $(a,b)\circ \phi=R(a)L(b)\phi$.

Recall \cite{D,R}
that the center of $U_\hh(\g)$ is spanned by elements 
$C_V$ corresponding to finite dimensional representations $V$ of $U_\hh(\g)$
by the formula 
$$
C_V=\text{Tr}|_V(1\o \pi_V)(\Cal R^{21}\Cal R(1\o e^{2\hh\rho})),\tag 3.1
$$
where $\Cal R$ is the universal R-matrix of $U_\hh(\g)$. 
The map $V\to C_V$ defines a homomorphism of the Grothendieck ring 
of the category of finite dimensional representations of $U_\hh(\g)$
and the center of $U_\hh(\g)$. 

Now 
we would like to define the notion of a Whittaker function. In order to do so, 
we will define an analogue of the notion of a nondegenerate character. 

Let $(a_{ij})$ be the Cartan matrix of $\g$. Let $d_i=1,2,3$ 
be the set of relatively prime integers such that the matrix 
$b_{ij}=d_ia_{ij}$ is symmetric. 
Let $\omega$ be any orientation of the Dynkin diagram of $\g$. 
Define the quantum polynomial algebra $P_\omega$ generated by 
variables $x_1,...,x_r$ corresponding to vertices of the Dynkin diagram, 
with relations $x_ix_j=e^{\pm \hh b_{ij}}x_jx_i$, where the sign is $+$ if
the edge $ij$ is oriented from $i$ to $j$ and $-$ otherwise. 

We will need the following well known proposition. 

\proclaim{Proposition 3.1} There exists a homomorphism 
$\chi_+:U_\hh(\n_+)\to P_\omega$ such that $\chi_+(e_i)=x_i$.
\endproclaim

\demo{Proof} It is sufficient to check that 
$\chi_+$ respects the Serre relations. Thus, it is enough to show that 
$$
\sum_{k=0}^{1-a_{ij}}(-1)^k\biggl[\matrix 1-a_{ij}\\ 
k\endmatrix\biggr]_{q^{d_i}}x_i^{1-a_{ij}-k}x_jx_i^k=0,\tag 3.2
$$
(here $q=e^\hh$, and $[a]_p=\frac{p^a-p^{-a}}{p-p^{-1}}$). 
Using the relations between $x_i$, we can reduce this to 
$$
\sum_{k=0}^{1-a_{ij}}(-1)^k\biggl[\matrix 1-a_{ij}\\
k\endmatrix\biggr]_{q^{d_i}}q^{\pm kb_{ij}}=0,
$$
which holds by the q-binomial theorem.
\enddemo

We also let $\chi_-:U_\hh(\n_-)\to P_\omega$ be the homomorphism 
defined by $\chi_-(f_i)=x_i$. 

Introduce the algebra $A_\omega=P_\omega\o A_{U_\hh(\g)}$. 
 
\proclaim{Definition} A Whittaker function on $U_\hh(\g)$ 
is an element $\phi\in A_\omega$ such that 
for any $a_\pm\in U_\hh(\n_\pm)$ one has 
$$
(a_-,a_+)\circ \phi=\chi_+(a_+)\phi\chi_-(a_-).\tag 3.3
$$
\endproclaim

It is obvious that the module $U_\hh(\g)$ over 
$U_\hh(\n_-)\o U_\hh(\n_+)$ (where 
the first factor acts by left shifts and the 
second by right shifts) is freely generated by the 
Cartan subalgebra $U_\hh(\h)$. Therefore,
the space of Whittaker functions is naturally identified  
with $P_\omega\o U_\hh(\h)$. As before, we will denote this identification 
by $\phi\to \phi|_\h$.

Let $Y\in U_\hh(\g)$. Then 
$Y$ defines an endomorphism $L(Y)$ of $A_{U_\hh(\g)}$ and 
$A_\omega$. As before, we have

\proclaim{Proposition 3.2} (i) For any $Y\in U_\hh(\g)$,
which is a noncommutative polynomial of elements $e^{\hbar x}$ ($x\in\h$), 
$e_i$, and $f_i$, there exists a 
unique difference operator $D_Y$ on $H$ with coefficients in 
$B=P_\omega\o P_\omega^{op}$ 
such that $(L(Y)\phi)|_\h=D_Y\phi|_\h$ for any Whittaker function $\phi$.

(ii) If $Y$ and $Y'$ are central 
then $D_{YY'}=D_YD_{Y'}$.
\endproclaim

\demo{Proof} Same as for Proposition 1.3. 
\enddemo

Now consider the fundamental representations 
$V_1,...,V_r$ of $U_\hh(\g)$ and define $Y_i=C_{V_i}$.
Then $D_{Y_i}$ 
are $B$-valued commuting difference operators 

Unfortunately, operators $D_{Y_i}$ are with 
coefficients in $B=P_\omega\o P_\omega^{op}$, 
while we want to obtain operators 
with scalar coefficients. However, since $C_V$ has zero weight for any $V$, 
the element $D_Y$ for $Y=C_V$ is in fact a difference operator with 
coefficients in the subalgebra $Q\subset P_\omega\o P_\omega^{op}$ generated 
by $x_i\o x_i$. It is easy to see that the algebra $Q$ is {\it commutative}. 
Thus, setting
$x_i\o x_i$ to be equal to any numbers $\beta_i\in\C$, we can 
obtain commuting difference operators $\tilde\Cal M_1,...,\tilde\Cal M_r$ 
with {\it scalar} coefficients. We will fix the normalization by 
letting $\beta_i=-1$, and
define an
system of $r$ commuting difference operators $\Cal M_i:=
e^{(\rho,h)}\tilde\Cal M_ie^{-(\rho,h)}$, $i=1,...,r$,

The fact that these difference operators are algebraically independent
follows from the facts that any central element of $U(\g)$ can be 
q-deformed, and that the differential operators of Section 1 are 
algebraically independent. 

We will call the system $\{\Cal M_i\}$ {\it the q-deformed Toda system}.
We note that this system depends on the choice of the orientation 
$\omega$ of the Dynkin diagram. 

\head{4. q-deformed affine Toda systems}\endhead

In this section we will generalize the constructions of Section 3 to 
quantum affine algebras. 

Consider the quantum affine algebra 
$U_\hh(\hat\g)$. 
It has generators $e_i$, $f_i$, $h_i$, $i=0,...,r$,  
and the standard relations \cite{CP}.
Let $U_\hh(\hat\n_+)$ be the subalgebra of $U_\hh(\hat\g)$ generated by $e_i$, 
and $U_\hh(\hat\n_-)$ be the subalgebra generated by $f_i$, respectively. 

Let $A_{U_\hh(\hat\g)}=U_\hh(\hat\g)^*$ be the 
space of linear functions on $U_\hh(\hat\g)$, 
i.e. the dual quantum formal group. 

The obvious two-sided action of $U_\hh(\hat\g)$ on itself 
induces a two-sided action of $U_\hh(\hat\g)$ on $A_{U_\hh(\hat\g)}$, via
$((a,b)\circ \phi)(x)=\phi(axb)$ (here the action of $a$ is a right 
action and the action of $b$ is a left action). Define 
the maps $L,R:U_\hh(\hat\g)\to \End(A_{U_\hh(\hat\g)})$ 
by $(a,b)\circ \phi=R(a)L(b)\phi$. 

For any complex number $k$, let $A^k_{U_\hh(\hat\g)}\subset A_{U_\hh(\hat\g)}$ 
be the space of functions 
satisfying the equation $L(c)f=kf$. 
Then it is clear that for any $Y\in U_\hh(\hat\g)$, 
$L(Y)$ preserves $A^k_{U_\hh(\hat\g)}$, and $L(c)|_{A^k_{U_\hh(\hat\g)}}=
k\cdot Id$. 
Thus $L$ descends to a map 
$L:U_\hh(\hat\g)/(c-k)\to \End_\C(A^k_{U_\hh(\hat\g)})$. 

As in the classical case, 
the value $k=-h^\vee$ is especially important,
since at this point there are a lot of interesting central 
elements. 
More specifically, one can generalize the Drinfeld-Reshetikhin 
construction of central elements (see Section 3) to the affine case. 
This was originally done in \cite{RS} (see also \cite{DE}). 
In this generalization, to every finite dimensional representation 
$V$ of $U_\hh(\hat\g)$ one assigns a central element 
$C_V\in \hat U_\hh(\hat\g)/(c+h^\vee)$, where 
$\hat U_\hh(\hat\g)$ is the completion of the quantum 
affine algebra acting in modules which are locally nilpotent 
under $e_i$. This defines a homomorphism 
from the Grothendieck algebra of the category of finite dimensional 
representations of $U_\hh(\hat\g)$ to $\hat U_\hh(\hat\g)$. 
The element $C_V$ is defined by the formula 
$$
C_V=Res_{z=0}Tr_V(Id\o \pi_{V(z)})(\Cal R^{21}\Cal R(1\o e^{2h\rho}))
\frac{dz}{z},\tag 4.1
$$
where $\Cal R$ is the truncated R-matrix of the quantum affine algebra 
(defined by formula (1.1) in \cite{DE}), and $V(z)$ is the representation $V$ 
shifted by $z\in \C^*$ (see \cite{DE}). 

Now 
we would like to define the notion of a Whittaker function. In order to do so, 
we will define an analogue of the notion of a nondegenerate character, 
like in Section 3.  

Let $(a_{ij})$ $(i,j\ge 0)$ be the Cartan matrix of $\hat\g$. Let $d_i=1,2,3$ 
be as in section 3 for $i>0$, and $d_0=1$. Let
$b_{ij}=d_ia_{ij}$. 
Let $\omega$ be an orientation of the Dynkin diagram of $\hat\g$.
Define the quantum polynomial algebra $P_\omega$ generated by 
variables $x_0,x_1,...,x_r$ corresponding to vertices of the Dynkin diagram, 
with relations $x_ix_j=e^{\pm \hh b_{ij}}x_jx_i$, where the sign is $+$ if
the edge $ij$ is oriented from $i$ to $j$ and $-$ otherwise. 

We will need the following proposition, which is the affine analog 
of Proposition 3.1. 

\proclaim{Proposition 4.1} There exists a homomorphism 
$\chi_+:U_\hh(\hat\n_+)\to P_\omega$ such that $\chi_+(e_i)=x_i$.
\endproclaim

\demo{Proof} Same as  
Proposition 3.1. $\square$\enddemo

We also let $\chi_-:U_\hh(\hat\n_-)\to P_\omega$ be the homomorphism 
defined by $\chi_-(f_i)=x_i$. 

Introduce the algebra $A_\omega=P_\omega\o A^{-h^\vee}_{U_\hh(\hat\g)}$. 
 
\proclaim{Definition} A Whittaker function on $U_\hh(\hat\g)$ 
is an element $\phi\in A_\omega$ such that 
for any $a_\pm\in U_\hh(\hat\n_\pm)$ one has 
$$
(a_-,a_+)\circ \phi=\chi_+(a_+)\phi\chi_-(a_-).\tag 4.2
$$
\endproclaim

As before, the space of Whittaker functions is naturally identified  
with $P_\omega\o U_\hh(\h)$. Namely, any element of $P_\omega\o U_\hh(\h)$ can 
be uniquely extended by equivariance to a Whittaker function. 
We denote this identification 
by $\phi\to\phi|_\h$. 

Now consider the fundamental representations 
$V_1,...,V_r$ of $U_\hh(\hat\g)$ (see \cite{CP}; 
they can be bigger than the fundamental representations 
of $U_\hh(\g)$ if $\g$ is not of type A), and define $Y_i=C_{V_i}$.

Let $\omega$ be an acyclic orrientation of the Dynkin diagram of $\hat\g$
(``acyclic'' is vacuous unless $\g$ is of type $A$).  

\proclaim{Proposition 4.2} 
(i) For any $i=1,...,r$ the series $(L(Y_i)\phi)|_\h$ is finite
for any Whittaker function $\phi$. Moreover, 
there exists  
a unique difference operator 
$D_{Y_i}:A_\h\to A_\h$ on $H$ such that for any Whittaker function 
$\phi$ one has $(L(Y_i)\phi)|_\h=D_{Y_i}\phi|_\h$. 

(ii) $[D_{Y_i},D_{Y_j}]=0$. 
\endproclaim

\demo{Proof} The finiteness of the series 
is obtained like in Section 2, using a suitable analog of 
the fact that terms containing non-simple roots 
vanish (see Section 5 for a detailed proof). 
The rest of the proof is analogous to the proof of Proposition 1.3. 
$\square$\enddemo

Thus $D_{Y_i}$ are 
$B$-valued commuting difference operators, where $B=P_\omega\otimes 
P_\omega^{op}$.  

As before, since $C_V$ has zero weight for any $V$, 
the element $D_Y$ for $Y=C_V$ is in fact a difference operator with 
coefficients in the subalgebra $Q\subset P_\omega\o P_\omega^{op}$ generated 
by $x_i\o x_i$. The algebra $Q$ is commutative. Thus, setting
$x_i\o x_i$ to be equal to any numbers $\beta_i\in\C$, we can 
obtain commuting difference operators $\tilde\Cal M^K_1,...,\tilde\Cal M^K_r$ 
with scalar coefficients. We will fix the normalization by 
letting $\beta_i=-1$, $i\ge 1$, $\beta_0=-K$, and 
define a
system of $r$ commuting difference operators $\Cal M^K_i=
e^{(\rho,h)}\tilde\Cal M^K_ie^{-(\rho,h)}$, $i=1,..,r$.

As before, the operators $\Cal M_i^K$ are algebraically independent. 
This follows from the fact that by Section 3 this is the case for $K=0$.

We will call the system $\{\Cal M^K_i\}$ 
{\it the q-deformed affine Toda system}.
It depends on the orientation 
on the Dynkin diagram.

\head 5. Computation of the q-deformed Toda operators\endhead

In this section we will complete the proof of Proposition 4.2, 
explain how to compute the q-deformed and q-deformed 
affine Toda operators, and compute some of them for $\g=sl(N)$. 

First of all, we need to recall the 
notion of Cartan-Weyl root elements $e_\beta$ of a finite dimensional 
or affine quantum group
(corresponding to all roots), 
due to Khoroshkin and Tolstoy \cite{KhT1,KhT2}.
To define them, one needs to fix a normal ordering of 
positive roots (see \cite{T}). 
(As was shown in \cite{T},
such an ordering can be obtained by extending any ordering of simple roots). 
Then one computes the corresponding elements $e_\beta$
by formulas (10-14) of \cite{KhT1}. 

\proclaim{Propsoition 5.1} Suppose that an orientation $\omega$ 
of the Dynkin diagram is 
consistent with the ordering of simple roots used 
to define the Cartan-Weyl root vectors (i.e. edges are oriented from 
smaller to larger simple roots). 
Then the homomorphisms $\chi_\pm$ corresponding to $\omega$ 
annihilate the Cartan-Weyl root vectors 
for non-simple roots.
\endproclaim

\demo{Proof}
This follows from the 
definition of the Cartan-Weyl generators, i.e. 
formulas (10-14)  in \cite{KhT1}. (The Cartan-Weyl generators are defined 
as iterated q-commutators of simple root elements, so they by definition 
map to $0$ under $\chi_\pm$). 
$\square$\enddemo

Next, we prove the following lemma about Whittaker functions.

We consider the quantum group $U_\hh(\g)$ or $U_\hh(\hat\g)$. 
Let us fix an acyclic orientation $\omega$ of its Dynkin diagram
and extend it to an ordering of simple positive roots. Let us further 
extend this ordering to a normal ordering of positive roots.  

\proclaim{Lemma 5.2} Let $X=f_{\gamma_1}...f_{\gamma_n}
e_{\gamma_1'}...e_{\gamma_m'}$, where 
$f_\gamma$, $e_\gamma$ are root elements from the Cartan-Weyl basis 
corresponding to roots $-\gamma$ and $\gamma$ ($\gamma>0$). If the roots
$\gamma_1,...,\gamma_n,\gamma_1',...,\gamma_m'$ are not all simple 
then $L(X)\phi=0$ for any Whittaker function $\phi$
(where to define Whittaker functions one used the orientation $\omega$).
Otherwise, if $\gamma_i=\alpha_{k_i}$, $\gamma_i'=\alpha_{l_i}$ then 
$$
(L(X)\phi)|_\h=e^{(\sum\alpha_{k_i},h)}x_{l_1}...x_{l_m}\phi|_\h
x_{k_n}...x_{k_1}\tag 5.1
$$ 
\endproclaim

\demo{Proof} By Proposition 5.1 we have: 
$L(e_\alpha)\phi$ equals $0$ if $\alpha$ is not simple 
and $x_i\phi$ if $\alpha=\alpha_i$. Therefore, it is easy to see 
that in order for $L(X)\phi$ to be nonzero, the roots 
$\gamma_i'$ must be simple. If they are (i.e. $\gamma_i'=\alpha_{l_i}$) then 
we have 
$$
(L(X)\phi)|_\h=x_{l_1}...x_{l_m}(L(f_{\gamma_1}...f_{\gamma_n})\phi)|_\h.
\tag 5.2
$$
Now, we have (like in the proof of Proposition 1.3): 
$$(L(f_{\gamma_1}...f_{\gamma_n})\phi)|_\h=
e^{(\sum \gamma_i,h)}(R(f_{\gamma_1}...f_{\gamma_n})\phi)|_\h.$$ 
Thus, similarly to the above, 
we get $0$ if at least one $\gamma_i$ is not simple. 
If $\gamma_i=\alpha_{k_i}$ then we get (5.1).
$\square$\enddemo

Now we are ready to compute the q-deformed Toda operators. 
The computation is based on the explicit formula 
for the R-matrix $\Cal R$ given in \cite{KhT1} (formula (42)). 
This formula says that 
$\Cal R$ (for finite dimensional or affine quantum groups) 
can be represented as a normally ordered product of factors corresponding 
to positive roots (for any fixed normal ordering). 

So let us fix a normal ordering as above and consider the 
Khoroshkin-Tolstoy representation of $\Cal R$. Let $\Cal R_s$ 
be the defined by the same product as $\Cal R$ but with 
all terms corresponding to non-simple roots crossed out. 
That is, 
$$
\Cal R_s=(\prod_i\Cal R_{\alpha_i})\Cal K,\tag 5.3
$$
where 
$\Cal K=q^{ \sum_{i=1}^ry_i\o y_i}$, for any 
orthonormal basis $y_i$ of $\h$, and $\Cal R_{\alpha_i}$ 
are R-matrices corresponding to simple roots:
$$
\Cal R_{\alpha_i}=
\exp_{q^{2d_i}}((q^{d_i}-q^{-d_i})e_{\alpha_i}\o f_{\alpha_i}),
\tag 5.4
$$
where 
$$
\exp_p(x)=\sum_{n\ge 0}\frac{(p-1)^n}{(p-1)...(p^n-1)}x^n\tag 5.5
$$
is the quantum exponential. 

{\bf Remark.}
The product in (5.3) is ordered according to the ordering 
of simple roots obtained by extension of $\omega$. It is not hard to show, 
however, that it depends only on $\omega$ itself and not on the extension
(This follows from the obvious fact that any two complete orders extending 
a partial order on a finite set can be identified by a sequence of 
transpositions of neighbors where transposed elements are not comparable).  

We have the following obvious corollary from Lemma 5.2.

\proclaim{Corollary 5.3} Let $C_V^s$ be defined by the same formula as $C_V$ 
(i.e. (3.1) or (4.1)) but with $\Cal R$ replaced by $\Cal R_s$. 
Let $Y_i^s=C_{V_i}^s$. Then $D_{Y_i^s}=D_{Y_i}$.  
\endproclaim

The operators $D_{Y_i^s}$ are relatively easy to compute in any given 
special case, since they contain
contributions from simple roots only. In the remainder of this section, 
we will restrict ourselves to $\g=sl(N)$ and let $V_i=\Lambda^iV$, where 
$V$ is the vector representation. We will compute 
the q-Toda operator corresponding to $i=1$ (i.e. $V_i=V$) explicitly.  

First of all, we will consider the non-affine case. 
We have $e_i^2=f_i^2=0$ in $V$, and so in the definition of $D_{Y_1^s}$ 
one can replace $\Cal R_{\alpha_i}$ with the first two terms
$1+X_i$, where $X_i=(q-q^{-1})e_i\o f_i$.  
Thus, $D_{Y_1^s}$ is obtained from 
a product of such binomial terms. After opening the brackets, 
we find that many terms are equal to zero. Indeed, 
we have 
$$
C_V=\sum_l\sum_{i_1,...,i_l}\text{Tr}|_V(X_{i_1}^{21}....X_{i_l}^{21}
X_{i_1}...X_{i_l}q^{2\rho})
$$
It is clear that a term in this sum can be nonzero 
only if $i_1,...,i_l$ are consequtive (as it contains $e_{i_1}...e_{i_l}|_V$)
and $i_l,...,i_1$ are consequtive (as it contains $f_{i_1}...f_{i_l}|_V$).
Thus, the terms with $l>1$ are zero. Computing the terms with $l=0,1$, 
we obtain
$$
C_V=\sum_{j=1}^Nq^{N+1-2j}q^{2\omega_j}+(q-q^{-1})^2
\sum_{i=1}^{N-1} q^{N+1-2i}f_iq^{\omega_{i+1}}
e_iq^{\omega_i},\tag 5.6
$$
where $\omega_i$ is the weight of the i-th basis vector of $V$. 
From (5.6), we get 
$$
\Cal M_1=\sum_{j=1}^NT_j^2-
(q-q^{-1})^2\sum_{i=1}^{N-1}
e^{(h,\alpha_i)}T_iT_{i+1},\tag 5.7
$$
where $T_i=T_{\omega_i}$. 

Now consider the affine case. In this case, using a similar argument to the 
above, instead of (5.6) we get 
$$
C_V=\sum_{j=1}^Nq^{N+1-2j}q^{2\omega_j}+(q-q^{-1})^2
\sum_{i=1}^{N} q^{N+1-2i}f_iq^{\omega_{i+1}}
e_iq^{\omega_i},\tag 5.8
$$
with subscripts understood modulo $N$, i.e. $0=N$. 
(here we use that, according to (4.1), we only take
the zero degree terms in the z-expansion of 
$(1\o \pi_{V(z)})(\Cal R^{21}\Cal R)$).  
Therefore, instead of (5.8) we get 
$$
\Cal M_1^K=\sum_{j=1}^NT_j^2-
(q-q^{-1})^2\sum_{i=1}^{N}K^{\delta_{iN}}
e^{(h,\alpha_i)}T_iT_{i+1}.\tag 5.9
$$

Similarly one can compute $\Cal M_i$, $\Cal M_i^K$, $i>1$
with a somewhat more complicated answer. 

\head 6. Quasiclassical limit of q-deformed quantum Toda systems, 
and their relation to quantum relativistic Toda systems\endhead

The following proposition 
 explains why we refer to the quantum integrable systems of Sections 3,4 
as q-deformed (non-affine or affine) Toda systems. 

\proclaim{Proposition 6.1} 
For any $i$ and $K$ one has
$$
\lim_{\hh\to 0}\frac{\Cal M^K_i-\text{dim}(V_i)}{(q-q^{-1})^2}=
C_i(M^K+G_i),\tag 6.1
$$
where $C_i\in\C^*$ and $G_i\in \C$.
\endproclaim

\demo{Proof} This proposition follows from Theorem 4.1 of \cite{DE}, 
which computes the quasiclassical limit of central elements. 
$\square$\enddemo

{\bf Remark.} In the special case $\g=sl(N)$, $i=1$, Proposition 6.1 can be 
checked directly from (5.9): it is easy to see that 
$$
\lim_{\hh\to 0}\frac{\Cal M^K_1-N}{(q-q^{-1})^2}=
-M^K_1,\tag 6.2
$$

Next, let us show that the q-deformed 
quantum Toda systems of Sections 3,4 (for $\g=sl(N)$) are equivalent to 
quantum relativistic Toda systems of \cite{Ru1}, nonperiodic and periodic, 
respectively. 

Introduce coordinates $z_1,...,z_N$ so that $\h$ is the set of solutions 
of $\sum z_i=0$. We will realize formal functions on $H$ as functions of 
$z_1,...,z_N$ invariant under simultaneous shift of $z_i$. 
In terms of $z_i$, the operator $\Cal M_1^K$ (as operator on such functions) 
can be written in the form
$$
\Cal M_1^K=\sum_{i=1}^NT_i^2-(q-q^{-1})^2
\sum_{i=1}^{N}K^{\delta_{iN}}e^{z_i-z_{i+1}}T_iT_{i+1},\tag 6.3
$$
where $T_if(z_1,...,z_i,...,z_N)=f(z_1,...,z_i+\hh,...,z_N)$
(we remind that subscripts are understood cyclically). 

Consider the algebra $Q$ of operators on
simultaneous-shift-invariant functions generated 
by $e^{\pm z_i\mp z_{i+1}}$ and $T_i^{\pm 1}$. It is clear that  
for all $i$ the operator $\Cal M_i^K$ belongs to $Q$. 
Consider the automorphism of $Q$ defined by 
$T_i\to T_i$, $e^{z_i-z_{i+1}}\to e^{z_i-z_{i+1}}T_iT_{i+1}^{-1}$. 
Under this automorphism, the operator $\Cal M_1^K$ is mapped to 
a simpler operator
$$
\bar\Cal M_1^K=
\sum_{i=1}^N\Cal T_i-(q-q^{-1})^2\sum_{i=1}^N
K^{\delta_{iN}}e^{z_i-z_{i+1}}\Cal T_i,\tag 6.4
$$
where $\Cal T_i:=T_i^2$.   

Thus, the q-deformed affine quantum Toda system for $sl(N)$ is 
equivalent (as an integrable system) to the system defined by the Hamiltonian 
(6.4). Setting $K=0$, we get that the q-deformed non-affine 
quantum Toda system is equivalent to the system defined by 
$$
\bar\Cal M_1=
\sum_{i=1}^N\Cal T_i-(q-q^{-1})^2\sum_{i=1}^{N-1}
e^{z_i-z_{i+1}}\Cal T_i,\tag 6.5
$$

Now recall from \cite{Ru1} 
that the Hamiltonian of the quantum relativistic Toda system 
is 
$$
\hat S_1=\sum_{i=1}^Nf(z_{i-1}-z_i)\Cal T_if(z_i-z_{i+1}),\tag 6.6
$$
where $f(a):=(1+g^2e^a)^{1/2}$, and $z_0=z_N$, $z_{N+1}=z_1$ in the periodic 
case, and $z_0=-\infty$, $z_{N+1}=+\infty$ in the non-periodic case. 
Let us conjugate $\hat S_1$ by the function 
$\prod_i\psi(z_i-z_{i+1})$, where $\psi$ satisfies the difference 
equation $\psi(x+2\hh)=\psi(x)f(x)^{-1}$, and the product is from $1$ to 
$N$ in the periodic case, and from $1$ to $N-1$ in the nonperiodic case.
  We get 
$$
\psi^{-1}\hat S_1\psi:=\bar S_1=\sum_{i=1}^N(1+g^2e^{z_i-z_{i+1}})T_i.\tag 6.7
$$
In the nonperiodic case, setting $g=\sqrt{-1}(q-q^{-1})$, we see that (6.7) 
becomes (6.5). In the periodic case, we set $g=\sqrt{-1}(q-q^{-1})K^{1/N}$, 
and $z_i'=z_i-\frac{i}{N}\ln K$, where $K\ne 0$ is arbitrary. 
Then (6.7) becomes (6.4). This demonstrates the equivalence of the 
q-deformed and relativistic Toda lattices.

{\bf Remark 1.} 
In the non-periodic relativistic Toda lattice, the parameter $g$ 
can be removed by a shift of variables, while in the periodic case 
it is an essential parameter. This corresponds to the presence 
of $K$ in the the affine and its absence in the non-affine case. 

{\bf Remark 2.} Note that our quantum group theoretic procedure 
of Sections 3 and 4 yields the Hamiltonian given by (6.3) and not 
the simpler Hamiltonian (6.4). In fact, one can get (6.4) instead of 
(6.3) by a slight modification of the procedure. 

Namely, if $H_i$ are any elements of $\h$ such that $\alpha_i(H_j)=
\alpha_j(H_i)$ then new elements $e_i'=e_ie^{\hh H_i}$ satisfy the 
quantum Serre relations and commutativity for orthogonal roots. 
Instead of the algebra $U_\hh(\hat\n_+)$ used in our argument, 
we could use $U_q(\hat\n_+')$, generated by $e_i'$, and similarly for 
$U_\hh(\hat\n_-)$. This would produce a different operator from (6.3), 
which is, however, equivalent to (6.3) by an automorphism 
of the algebra of difference operators. In 
particular, it is easy to see that one can choose $H_i$ so that 
the obtained operator is (6.4). 

In this connection we would like to mention the paper \cite{S}, where 
a similar idea is used: 
the author finds elements $H_i$ such that the elements $e_i'$ 
generate an algebra which has nondegenerate characters 
(such $H_i$ are determined from a system of linear nonhomogeneous 
equations, whose homogeneous part is the above system 
$\alpha_i(H_j)=\alpha_j(H_i)$). Using such characters, 
one may also define q-deformed quantum Toda systems.
This approach is closely related to ours. For example, in 
\cite{S} the system for $H_i$ depends on a choice 
of a Coxeter element, which in our situation corresponds to
the necessity to choose an orientation of the Dynkin diagram.  

\head 7. Toda systems as limits of Calogero-Moser, 
Macdonald, and Ruijsenaars systems.\endhead

In this section we discuss the limiting procedures which allow 
to obtain quantum Toda systems as limiting cases of more complicated 
integrable systems inbvolving all roots of a Lie algebra rather than just 
simple roots. The results presented here are known and are given only for 
the sake of completeness of the picture. 

It is known \cite{I} that quantum Toda systems 
can be represented as limits of quantum Calogero-Moser systems.

Namely, for any finite dimensional Lie algebra $\g$ consider 
the quantum trigonometric Calogero-Moser Hamiltonian
$$
H_T(k)=-\frac{1}{2}\Delta+k(k-1)
\sum_{\alpha>0}\frac{1}{\text{sinh}^2\alpha(h)}.
\tag 7.1
$$
Set $k=\frac{1}{2}e^{P}$, $h=-\frac{1}{2}x+P\rho$. In terms of 
the new notation, we have
$$
H_T(k)=-\frac{1}{2}\Delta+e^{P}(e^P-2)\sum_{\alpha>0}
\frac{1}{4\text{sinh}^2(-\frac{1}{2}\alpha(x)+P(\alpha,\rho))},
\tag 7.2
$$
This implies that 
$$
\lim_{P\to+\infty}H_T(k)=
-\frac{1}{2}\Delta+\sum_{i=1}^re^{\alpha_i(h)}=M.
\tag 7.3
$$
(Indeed, only terms corresponding to simple roots remain finite:
the term corresponding to a root $\alpha>0$ behaves like 
$\text{const}\cdot e^{2P(1-(\alpha,\rho))}$). 
This proves our statement in the non-affine case. 

For the affine case, we consider the elliptic Calogero-Moser Hamiltonian 
$$
H_E(k)=-\frac{1}{2}\Delta-4\pi^2 k(k-1)\sum_{\alpha>0}
\wp(\frac{\alpha(h)}{2\pi i},it),
\tag 7.4
$$
where $\wp$ is the Weierstrass elliptic function with periods $1,it$, and 
$t>0$. Set 
$k=\frac{1}{2}e^{P}$, and 
$h=-\frac{1}{2}x+P\rho$, $t=Ph^\vee-\frac{1}{2}\ln K$, where $h^\vee$ 
is the dual Coxeter number of $\g$. Recall that 
$$
-4\pi^2 \wp(\frac{z}{2\pi i},it)=\sum_{n\in\Z}\text{sinh}^{-2}(z+tn)
+\gamma(t).\tag 7.5
$$ 
Let $m$ be the number of positive roots of $\g$. 

We have 
$$
\gather
H_E(k)-mk(k-1)\gamma(t)=\\
-\frac{1}{2}\Delta+\frac{1}{4}e^{P}(e^P-2)\sum_{\alpha>0,n\in\Z}
\text{sinh}^{-2}(-\frac{1}{2}\alpha(x)+P(\alpha,\rho)+nPh^\vee
-\frac{n}{2}\ln K)\tag 7.6\endgather
$$
(for simplicity we assume that $K>0$). 

Therefore, 
$$
\lim_{P\to +\infty}(H_E(k)-mk(k-1)\gamma(t))=
-\frac{1}{2}\Delta+\sum_{i=0}^rK^{\delta_{i0}}e^{\alpha_i(h)}=M^K.\tag 7.7
$$
(Now the surviving roots are not only the simple roots of $\g$ but also the 
maximal root, because by the definition  $h^\vee=1+\theta(\rho)$). 
This proves our claim in the affine case.

This limiting procedure can be used to give another proof of Theorem 2.1 
which does not use quantum groups. Namely, it was shown by Cherednik \cite{Ch1}
that the quantum system defined by the Hamiltonian (7.4) is integrable 
for all Lie algebras. By a limiting argument, one can deduce 
from this that the same is true for the Hamiltonian $M^K$. 

Similarly, it is known that q-deformed quantum Toda systems
for $\g=sl(N)$ can be viewed as 
limits of quantum Macdonald-Ruijsenaars 
systems defined in \cite{Mac,Ru2}. 

Consider the non-affine case. 
Recall \cite{Mac,Ru2} that the trigonometric Macdonald-Ruijsenaars system
is defined by the quantum Hamiltonian 
$$
H_T^q=\sum_{i=1}^N(\prod_{j\ne i}\frac{q^{2k}e^{w_i}-e^{w_j}}
{e^{w_i}-e^{w_j}})\Cal T_i.
\tag 7.8
$$
Let us conjugate $H_T^q$ with the function 
$\eta=e^{-P\sum_{i=1}^N(i-1)w_i}$.
We get 
$$
\hat H_T^q=\eta^{-1}H_T^q\eta=
\sum_{i=1}^Nq^{-2(i-1)P}(\prod_{j\ne i}\frac{q^{2k}e^{w_i}-e^{w_j}}
{e^{w_i}-e^{w_j}})\Cal T_i.\tag 7.9 
$$
Set $w_i=z_i+2\hh i P$, $k=P$. We get 
$$
\hat H_T^q=
\sum_{i=1}^Ne^{-2\hh (i-1) P}(\prod_{j\ne i}
\frac{e^{2\hh (i+1)P}e^{z_i}-e^{2\hh jP}
e^{z_j}}
{e^{2\hh iP}e^{z_i}-e^{2\hh jP}e^{z_j}})\Cal T_i.\tag 7.10 
$$
Now let $\hh P\to +\infty$ (here $\hh$ has to be a complex number, 
not a formal 
parameter). It is easy to see that  
$$
\hat H_T^q\to \Cal T_N+\sum_{i=1}^{N-1}(1-e^{z_i-z_{i+1}})\Cal T_i.\tag 7.11
$$
This coincides with (6.5) after a shift of variables $z_i$. 

A similar computation shows that the Hamiltonian given by (6.4) 
can be obtained by a limiting procedure from the Ruijsenaars'
relativistic elliptic Calogero-Moser Hamiltonian \cite{Ru2}.

{\bf Remark 1.}
We expect that similar results are the case for all Lie algebras $\g$.
Namely, we expect that the non-affine and affine q-Toda systems
can be obtained as a limit of trigonometric, respectively elliptic 
Macdonald-Ruijsenaars operators, which were defined 
for an arbitrary root system by Cherednik \cite{Ch2,Ch3}.

{\bf Remark 2.} 
In \cite{E}, it is shown that Calogero-Moser operators 
(trigonometric and elliptic) are obtained as 
radial parts of central elements of $U(\g)$ (respectively, of a 
completion of $U(\hat\g)$ at the critical level) 
on equivariant functions on $G$ (respectively, $\hat G$)
 with values in certain 
special representations $U_k$. 
In \cite{EK}, it is shown that Macdonald operators 
are obtained in a similar manner 
as radial parts of central elements of $U_\hh(\g)$. 
In view of these results and the results of this paper, 
it would be tempting to understand
the above limiting procedures in terms of representation 
theory (i.e. to see how equivariant functions with values in $U_k$ 
turn into Whittaker functions in the limit).

\end